\documentclass[10pt, reqno]{amsart}
\usepackage[utf8]{inputenc}
\usepackage{amsmath}
\usepackage{hyperref}
\usepackage{amsfonts}
\hypersetup{
    colorlinks=false,
    linkcolor=blue,
    }

\usepackage{latexsym}
\usepackage{fancyhdr,amssymb}

\usepackage{amsrefs}
\usepackage{graphicx}
\graphicspath{ {./images/} }

\title{Maclaurin Integration:\\ A weapon against infamous integrals}
\author{Glenn Bruda}
\date{January 2022}

\DeclareMathOperator{\si}{Si}
\DeclareMathOperator{\ei}{Ei}

\usepackage{mathtools}

\DeclarePairedDelimiter\abs{\lvert}{\rvert}%
\DeclarePairedDelimiter\norm{\lVert}{\rVert}%

\makeatletter
\let\oldabs\abs
\def\abs{\@ifstar{\oldabs}{\oldabs*}}
\let\oldnorm\norm
\def\norm{\@ifstar{\oldnorm}{\oldnorm*}}
\makeatother

\begin{document}

\fancyhead{}
\fancyfoot{}
\fancyfoot[LE,RO]{\medskip \thepage}

\setcounter{page}{1}

\thanks{Special thanks to John Pfeilsticker and Kevin Knudson for your guidance and suggestions throughout the publication process.}

\maketitle

\begin{flushleft}
Perhaps one of the most challenging aspects of integration in calculus is that there is not a universal integration technique that works for all integrals. This paper introduces a technique by series to help solve this problem. \linebreak

To introduce the formula (\ref{formula}) and the nature of it, we consider two integrals and how they are approached classically:

\begin{enumerate}
    \item 
        \begin{align*}
            \int\frac{\sin(x)}{x}dx ~\text{and}
        \end{align*}
    \item     
        \begin{align*}
            \int e^{e^{x}}dx.
        \end{align*}

\end{enumerate}

Traditional integration techniques you may have learned in your calculus courses will not be able to solve either of these. These integrals can be calculated using the non-elementary functions $\si(x)$~\cite{trigint} and $\ei(x)$~\cite{expint} respectively, or by using a power series.
\vspace{\baselineskip}

A great strength of this formula, which we shall introduce, is that in addition to being able to solve infamously difficult integrals (like the ones above), is its ease of use. Compared to other integration techniques such as Trigonometric substitution, Integration by Partial Fractions, or Integration by Parts, Maclaurin Integration requires by far the least amount of labor to utilize. All that is needed to solve an integral using this technique is plugging in $f(x)$ for the integrand.
\vspace{\baselineskip}

\fancyhead{}
\fancyhead[CO]{\hfill MACLAURIN INTEGRATION}
\fancyhead[CE]{G. BRUDA  \hfill}
\renewcommand{\headrulewidth}{0pt}

To demonstrate this, we consider the integral
\begin{align*}
    \int e^{x^2}dx.
\end{align*}

The most common approach to evaluating this integral is to expand it as a power series and integrate term-by-term, which yields
\begin{align*}
    C + x +\frac{x^3}{3}+\frac{x^5}{10}+\frac{x^7}{42}+\frac{x^9}{216}+\cdots
\end{align*}
as the antiderivative, with $C$ as the constant of integration.
\linebreak
Maclaurin Integration is an alternative solution by series (although not a power series, since it involves the function itself in the solution) that eliminates the nuisance of calculation completely. To evaluate the integral above, simply substitute $e^{x^2}$ for $f(x)$ in the formula (see Eq. (\ref{formula})) to calculate the integral (see Eq. (\ref{3.2.1})).

\end{flushleft}

\centering
\section{Proof}

\subsection{Representing a Given Function with a Series}

\begin{flushleft}
To begin the proof of this formula, start with the Maclaurin Series
\begin{equation}
    \frac{1}{1-x} = \sum_{n=0}^{\infty} x^{n}.\tag{1.1.1}
\end{equation}
    Note the domain as $(-1, 1)$. Substituting $1-x$ for $x$,
\begin{equation}
    \frac{1}{x} = \sum_{n=0}^{\infty} (1-x)^{n}\notag~\text{and}
\end{equation}
\begin{equation}
    1 = x\sum_{n=0}^{\infty} (1-x)^{n}\notag.
\end{equation}
    Note the domain now as $(0, 2)$. Multiplying by $f(x)$,
\begin{equation}
    f(x) = xf(x)\sum_{n=0}^{\infty} (1-x)^{n}\notag~\text{and}
\end{equation}
\begin{equation}
    \int{f(x)}dx = \int{xf(x)\sum_{n=0}^{\infty} (1-x)^{n}}dx\tag{1.1.2}\label{maclaurinfx}.
\end{equation}
\end{flushleft}

\subsection{Applying Integration by Parts}
\pagebreak
\begin{flushleft}
Given the Integration by Parts Formula:
\begin{equation}
    \int{u}{dv} = uv-\int{vdu},\tag{1.2.1}
\end{equation}
substitute $xf(x)$ for $u$ and $\sum_{n=0}^{\infty} (1-x)^{n}$ for $dv$. We now apply Integration by Parts to Eq. (\ref{maclaurinfx}).
\begin{align*}
    \int{xf(x)\sum_{n=0}^{\infty}(1-x)^{n}}dx = xf(x)\sum_{n=0}^{\infty} \frac{-(1-x)^{n+1}}{n+1}    \\  -\int{\sum_{n=0}^{\infty} \frac{-(1-x)^{n+1}}{n+1}\left(\frac{d}{dx}(xf(x))\right)dx}\notag.
\end{align*}
Hence,
\begin{equation}
    \int{f(x)dx} = xf(x)\sum_{n=0}^{\infty} \frac{-(1-x)^{n+1}}{n+1}     -\int{\sum_{n=0}^{\infty} \frac{-(1-x)^{n+1}}{n+1}\left(\frac{d}{dx}\left(xf(x)\right)\right)dx}\notag
\end{equation}
\begin{equation}
    = -xf(x)\sum_{n=0}^{\infty} \frac{(1-x)^{n+1}}{n+1}+\int\sum_{n=0}^{\infty} \frac{(1-x)^{n+1}}{n+1}\left(\frac{d}{dx}\left(xf(x)\right)\right)dx\notag.
\end{equation}
We will continue applying Integration by Parts to the tail integral until we establish a clear pattern.
\linebreak
\linebreak
Substituting $\frac{d}{dx}(xf(x))$ for $u$ and $\sum_{n=0}^{\infty} \frac{(1-x)^{n+1}}{n+1}$ for $dv$ yields
\begin{align*}
   \int{f(x)dx} = -xf(x)\sum_{n=0}^{\infty} \frac{(1-x)^{n+1}}{n+1} - \frac{d}{dx}(xf(x))\sum_{n=0}^{\infty} \frac{(1-x)^{n+2}}{(n+1)(n+2)} \\ + \int\frac{d^2}{dx^2}(xf(x))\sum_{n=0}^{\infty}\notag \frac{(1-x)^{n+2}}{(n+1)(n+2)}dx.
\end{align*}
Substituting $\frac{d^2}{dx^2}(xf(x))$ for $u$ and $\sum_{n=0}^{\infty} \frac{(1-x)^{n+2}}{(n+1)(n+2)}$ for $dv$ yields
\begin{align*}
   \int{f(x)dx} = -xf(x)\sum_{n=0}^{\infty} \frac{(1-x)^{n+1}}{n+1} - \frac{d}{dx}(xf(x))\sum_{n=0}^{\infty} \frac{(1-x)^{n+2}}{(n+1)(n+2)}   \\ - \frac{d^2}{dx^2}(xf(x))\sum_{n=0}^{\infty} \frac{(1-x)^{n+3}}{(n+1)(n+2)(n+3)} \\ + \int\frac{d^3}{dx^3}(xf(x))\sum_{n=0}^{\infty} \frac{(1-x)^{n+3}}{(n+1)(n+2)(n+3)}dx\tag{1.2.2} \label{4termsum}.
\end{align*}
From here, we can see a clear pattern has been established, and can move on to representing this sum with a series.
\end{flushleft}

\subsection{Eliminating the Tail Term}

\begin{flushleft}

We now demonstrate the Integration by Parts pattern with the series

\begin{align*}
    \int f(x)dx = \lim_{N\to\infty}\Biggl(-\sum_{u=0}^{N}\left(\frac{d^u}{dx^u}(xf(x))\right)\sum_{n=0}^{\infty}\frac{(1-x)^{n+u+1}}{\prod_{v=1}^{u+1}(n+v)} \\ + \int\frac{d^{N+1}}{dx^{N+1}}(xf(x))\sum_{n=0}^{\infty}\frac{(1-x)^{n+N+2}}{\prod_{v=1}^{N+2}(n+v)}dx + C\Biggl) \notag .
\end{align*}

It seems that there has been no progress eliminating the integral sign, but in fact the tail term with the integral can be completely omitted from the equation. This is because the last term is zero. \linebreak

To prove this, we dissect the tail term. In many circumstances, the limit $\lim_{N\to\infty}\left(\frac{d^N}{dx^{N}}(xf(x))\right)$ will be effectively zero, but this is not always the case (ex: $f(x)=\frac{e^x}{x}$, where that limit yields $e^x$). What we need to evaluate is whether
\begin{align*}
    \lim_{N\to\infty}\left(\sum_{n=0}^{\infty}\frac{(1-x)^{n+N+2}}{\prod_{v=1}^{N+2}(n+v)}\right) = 0 \label{1.3.1} \tag{1.3.1}
\end{align*}

in the domain of $(0, 2)$. Factoring Eq. (\ref{1.3.1}) yields
\begin{align*}
    \lim_{N\to\infty}\left(\left(\sum_{n=0}^{\infty}\frac{(1-x)^{n+2}}{\prod_{v=1}^{N+2}(n+v)}\right)(1-x)^{N}\right).
\end{align*}
We separately evaluate the limit $\lim_{N\to\infty}((1-x)^N)$.
If $\lim_{N\to\infty}(x^N) = 0$ in the domain of $(-1, 1)$~\cite{limit}, then $\lim_{N\to\infty}((1-x)^N) = 0$ in the domain of $(0, 2)$ (as substituting $1-x$ for $x$ will change the applicable domain). Given this, this means that Eq. (\ref{1.3.1}) is true if \begin{align*}
    \lim_{N\to\infty}\left(\sum_{n=0}^{\infty}\frac{(1-x)^{n+2}}{\prod_{v=1}^{N+2}(n+v)}\right)
\end{align*}
is absolutely convergent. To determine this, we use the Ratio Test.

\begin{align*}
    \lim_{n\to\infty}\abs{\lim_{N\to\infty}\left(\left(\frac{(1-x)^{n+3}}{\prod_{v=1}^{N+2}(n+v+1)}\right)\cdot\left(\frac{\prod_{v=1}^{N+2}(n+v)}{(1-x)^{n+2}}\right)\right)}=L
\end{align*}

\begin{align*}
    =\lim_{n\to\infty}\abs{\lim_{N\to\infty}\left((1-x)\left(\prod_{v=1}^{N+2}\frac{n+v}{n+v+1}\right)\right)}.
\end{align*}
In order to continue to simplify our above limit, we individually evaluate the limit of the product operator:
\begin{align*}
    \lim_{n\to\infty}\abs{\lim_{N\to\infty}\left(\prod_{v=1}^{N+2}\frac{n+v}{n+v+1}\right)}. \tag{1.3.2} \label{1.3.2}
\end{align*}

Evaluating this, we find it to be a telescoping product: 
\begin{align*}
    \prod_{v=1}^{N+2}\frac{n+v}{n+v+1}=\left(\frac{n+1}{n+2}\right)\left(\frac{n+2}{n+3}\right)\left(\frac{n+3}{n+4}\right)\left(\frac{n+4}{n+5}\right)\cdots=\frac{n+1}{n+N+3}.
\end{align*}

Substituting the telescoping product for the product in Eq. (\ref{1.3.2}) yields
\begin{align*}
    \lim_{n\to\infty}\abs{\lim_{N\to\infty}\left(\frac{n+1}{n+N+3}\right)}.
\end{align*}
\begin{align*}
    \lim_{N\to\infty}\left(\frac{n+1}{n+N+3}\right) = 0 ~\text{and}
\end{align*}
\begin{align*}
     \lim_{N\to\infty}\left(\prod_{v=1}^{N+2}\frac{n+v}{n+v+1}\right) = 0, ~\text{which assures that}
\end{align*}
\begin{align*}
    \lim_{n\to\infty}\abs{\lim_{N\to\infty}\left((1-x)\left(\prod_{v=1}^{N+2}\frac{n+v}{n+v+1}\right)\right)} = 0.
\end{align*}
Since $L<1$, we can say that 
\begin{align*}
    \lim_{N\to\infty}\left(\sum_{n=0}^{\infty}\frac{(1-x)^{n+2}}{\prod_{v=1}^{N+2}(n+v)}\right)
\end{align*}
is absolutely convergent, thus proving Eq. (\ref{1.3.1}) true. This completes the proof that the tail term is 0.

\end{flushleft}

\section{Formula} \label{formula}

\begin{flushleft}

Since the tail term is omitted, the formula is simply:

\begin{align*}
    \int f(x)dx = -\sum_{u=0}^{\infty}\left(\frac{d^u}{dx^u}(xf(x))\sum_{n=0}^{\infty}\frac{(1-x)^{n+u+1}}{\prod_{v=1}^{u+1}(n+v)}\right) + C, \notag
\end{align*}

where $C$ is the constant of integration.

\end{flushleft}

\subsection{Conditions} \label{conditions}

\begin{flushleft}

This formula is valid only if:

\begin{enumerate}
    \item $f(x)$ is defined on the domain $(0, 2)$,
    \item $f(x)$ is continuous on (0, 2), and
    \item $xf(x)$ has derivatives of all orders on $(0,2)$.
\end{enumerate}

\end{flushleft}

\section{Applications}
\subsection{The Gamma Function}

\begin{flushleft}

Since this set of conditions are fairly liberal, this formula will work for many integrands. My initial idea for an application for this formula was the Gamma function (extension of $x!$ to complex numbers)~\cite{gamma}.
\linebreak \linebreak

The Gamma function $\Gamma(z)$ is equal to $(z-1)!$ when $z$ is any positive integer greater than 0. \linebreak

The most common approach to integrating the Gamma function is simply by integrating the integral that defines the Gamma function itself. There is nothing wrong with this approach, but it has its drawbacks. First, the answer is simply another integral. Further, the process to attain the indefinite integral of the Gamma function is not particularly difficult, but it does require more labor than integrating by Maclaurin Integration.
\linebreak

Note that it is possible to plug $\Gamma(x)$ in for $f(x)$ in the formula, but here we will present an alternative which does not use non-elementary functions.
\linebreak

First, we define the Gamma function using elementary functions. We will use Weierstrass's definition~\cite{weier}
\begin{align*}
    \Gamma(z) = \frac{e^{-\gamma z}}{z}\prod_{n=1}^{\infty}\left(\left(1+\frac{z}{n}\right)^{-1}e^{\frac{z}{n}}\right),
\end{align*}

where $\gamma$ is the Euler-Mascheroni constant~\cite{mascheroni}. Next we substitute Weierstrass's definition for $f(x)$ into the formula, while replacing $z$ with $x$. The resulting integral of $\Gamma(x)$ is

\begin{align*}
    -\sum_{u=0}^{\infty}\left(\frac{d^u}{dx^u}\left(e^{-\gamma x}\prod_{n=1}^{\infty}\left(\left(1+\frac{x}{n}\right)^{-1}e^{\frac{x}{n}}\right)\right)\sum_{n=0}^{\infty}\frac{(1-x)^{n+u+1}}{\prod_{v=1}^{u+1}(n+v)}\right) + C.
\end{align*}

\end{flushleft}

\subsection{Special Integrals}

\begin{flushleft}

For this subsection, we will consider four integrals whose antiderivatives cannot be represented using elementary functions. Integrals 1 and 2 are found using the error function~\cite{error}, while integrals 3 and 4 are found using the Fresnel integrals~\cite{fresnel}.

\begin{enumerate}
    \item 
        \begin{align*}
            \int e^{x^{2}}dx,
        \end{align*}
    \item     
        \begin{align*}
            \int e^{-x^{2}}dx,
        \end{align*}
    \item     
        \begin{align*}
            \int \sin({x^2})dx, ~\text{and}
        \end{align*}
    \item     
        \begin{align*}
            \int \cos({x^2})dx.
        \end{align*}
\end{enumerate}

Since all of these integrands satisfy the conditions (\ref{conditions}) for the formula to be valid, we can integrate all of these using the formula.

\begin{enumerate}
    \item 
        \begin{align*}
            \int e^{x^{2}}dx = -\sum_{u=0}^{\infty}\left(\frac{d^u}{dx^u}(xe^{x^{2}})\sum_{n=0}^{\infty}\frac{(1-x)^{n+u+1}}{\prod_{v=1}^{u+1}(n+v)}\right) + C. \label{3.2.1} \tag{3.2.1} 
        \end{align*}
    \item     
        \begin{align*}
            \int e^{-x^{2}}dx = -\sum_{u=0}^{\infty}\left(\frac{d^u}{dx^u}(xe^{-x^{2}})\sum_{n=0}^{\infty}\frac{(1-x)^{n+u+1}}{\prod_{v=1}^{u+1}(n+v)}\right) + C.
        \end{align*}
    \item     
        \begin{align*}
            \int \sin({x^2})dx = -\sum_{u=0}^{\infty}\left(\frac{d^u}{dx^u}(x\sin({x^2}))\sum_{n=0}^{\infty}\frac{(1-x)^{n+u+1}}{\prod_{v=1}^{u+1}(n+v)}\right) + C.
        \end{align*}
    \item     
        \begin{align*}
            \int \cos({x^2})dx = -\sum_{u=0}^{\infty}\left(\frac{d^u}{dx^u}(x\cos({x^2}))\sum_{n=0}^{\infty}\frac{(1-x)^{n+u+1}}{\prod_{v=1}^{u+1}(n+v)}\right) + C.
        \end{align*}
\end{enumerate}

\end{flushleft}

\subsection{Graphical Behavior}

\begin{flushleft}

In this subsection, we present an example function to find an antiderivative of using the formula. We will also render this result graphically in order to demonstrate the behavior of the formula at different upper limits of summation. \linebreak For our example, we will graph the antiderivative of
\begin{align*}
    \frac{x^5}{x^7+1} \tag{3.3.1} \label{3.3.1}
\end{align*}
using Maclaurin Integration. Normally, this an extremely difficult antiderivative to graph, but using this technique, it can be graphed accurately without an egregious amount of computation (nor large amount of human effort).

\includegraphics[width=10cm, height=10cm]{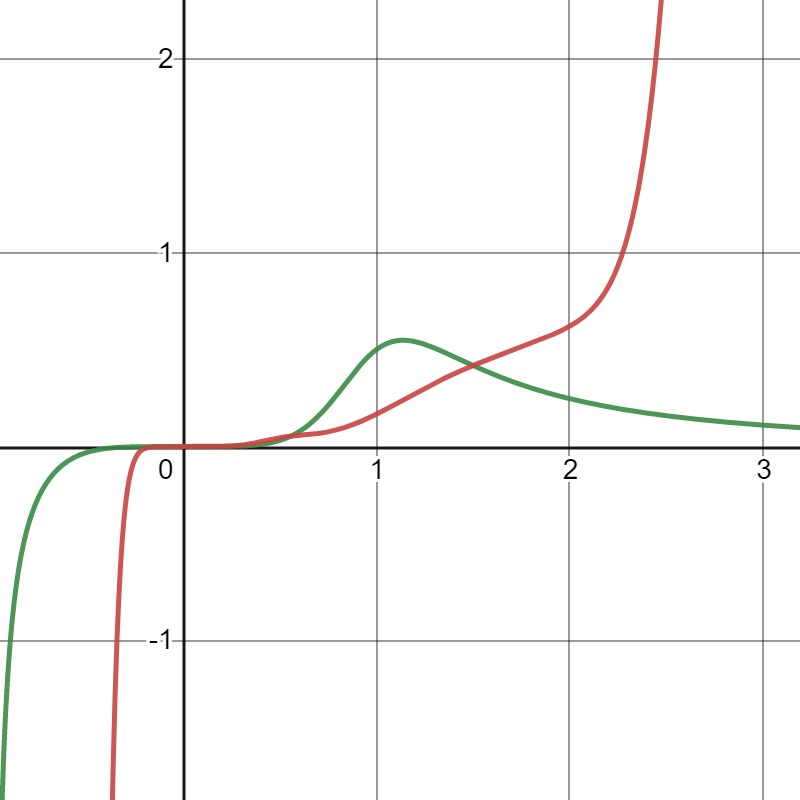}

The graph of the function $\frac{x^5}{x^7+1}$ is shown in green, and the graph of
\begin{align*}
    -\sum_{u=0}^{6}\left(\frac{d^u}{dx^u}\left(\frac{x^6}{x^7+1}\right)\sum_{n=0}^{10}\frac{(1-x)^{n+u+1}}{\prod_{v=1}^{u+1}(n+v)}\right)
\end{align*}
is shown in red.

Note that each infinite series' upper limits of summation in the formula were altered to a finite number for rendering and calculation purposes. \linebreak 

We now consider a simple definite integral of the integrand \ref{3.3.1} within the domain $(0,2)$, namely:
\begin{align*}
    \int_{1/2}^{3/2}\frac{x^5}{x^7+1}dx.
\end{align*}
We find this value to be approximately $0.3698$~\cite{wolframalpha}. 

To demonstrate the accuracy of the formula, we will make a table, but first consider the equation
\begin{align*}
    M(x)=-\sum_{u=0}^{p}\left(\frac{d^u}{dx^u}\left(\frac{x^6}{x^7+1}\right)\sum_{n=0}^{10}\frac{(1-x)^{n+u+1}}{\prod_{v=1}^{u+1}(n+v)}\right)dx.
\end{align*}
In our aforementioned table's first column, we will represent the upper limit of summation for $u$ (represented by $p$), and in the second column we will show the value of $M(\frac{3}{2})-M(\frac{1}{2})$ rounded to four decimal places. This should approach the value of the definite integral $\int_{1/2}^{3/2}\frac{x^5}{x^7+1}dx$ as $p$ grows larger.
\begin{center}
\begin{tabular}{|c| c|} 
 \hline
 Upper Limit of $u$ ($p$) & $M(\frac{3}{2})-M(\frac{1}{2})$ \\ [0.5ex] 
 \hline\hline
 0 & 0.2661 \\ 
 \hline
 1 & 0.3222 \\
 \hline
 2 & 0.3575 \\
 \hline
 3 & 0.3852 \\
 \hline
  4 & 0.3942 \\
 \hline
  5 & 0.3862 \\
 \hline
 6 & 0.3733 \\ [1ex] 
 \hline
\end{tabular}
\end{center}

As the table shows, the function $M(x)$ appears to become more accurate as $p$ increases, but not in a straightforward fashion. Instead, the values oscillate around the known approximate value of $0.3698$, as shown by the following graph (with $p$ on the x-axis, and the value of $M(\frac{3}{2})-M(\frac{1}{2})$ on the y-axis):
\includegraphics[width=12cm, height=5.5cm]{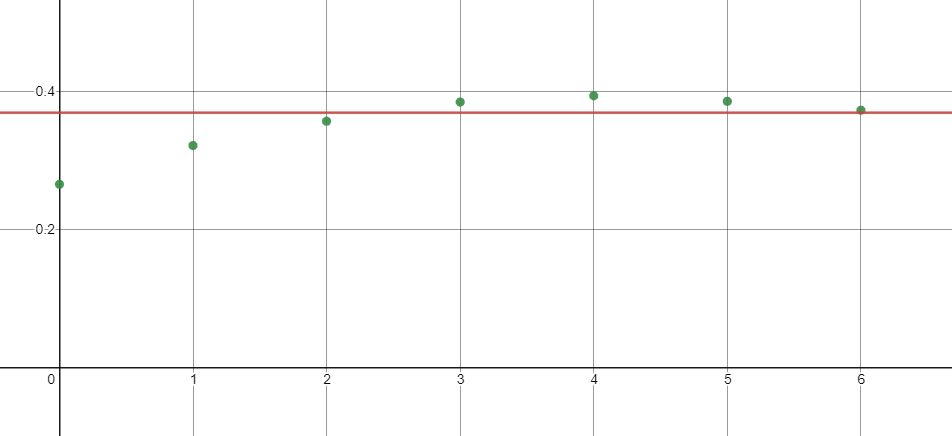}

with the red line indicating the value of $\int_{1/2}^{3/2}\frac{x^5}{x^7+1}dx$, and the green points being from the above table. While each point is not necessarily more accurate than the last (ex: at $p=4$, $M(x)$ is further from the accurate value than $M(x)$ is at $p=3$), $M(x)$ does become more accurate over larger periods of $p$ (ex: from $p=1$ to $p=6$). This demonstrates that, as $p$ grows larger, $M(x)$ does approach the value of the definite integral $\int_{1/2}^{3/2}\frac{x^5}{x^7+1}dx$, implying that the process is accurate. Note that this accuracy is achieved while having upper limits of summation for $n$ and $u$ equal to 10 or less. Graphical accuracy will greatly improve if the upper limits of summation are increased.

\end{flushleft}

\pagebreak

\begin{flushleft}
\vspace{\baselineskip}
\end{flushleft}


\begin{bibdiv}
\begin{biblist}

\bib{trigint}{webpage}{
    accessdate={2021-8-20},
    date={2015-06-11},
    title={Sine and Cosine Integrals},
    url={http://www.netlib.org/math/docpdf/ch02-14.pdf},
}

\bib{expint}{webpage}{
    accessdate={2021-8-20},
    date={2015-06-11},
    title={Exponential Integrals Ei and $E_{1}$},
    url={http://www.netlib.org/math/docpdf/ch02-10.pdf},
}

\bib{limit}{webpage}{
    accessdate={2021-8-20},
    title={The Limit of a Sequence},
    url={https://math.mit.edu/~apm/ch03.pdf},
    note={Theorem 3.4}
}

\bib{gamma}{webpage}{
    accessdate={2021-8-5},
    author={Sebah, Pascal},
    author={Gourdon, Xavier},
    date={2002-2-4},
    title={Introduction to the Gamma Function},
    url={http://scipp.ucsc.edu/~haber/archives/physics116A10/gamma.pdf},
    note={Theorem 2}
}

\bib{weier}{webpage}{
    accessdate={2021-8-22},
    title={Gamma Function},
    url={https://www.math.lsu.edu/system/files/WM1\%20paper.pdf},
}

\bib{mascheroni}{webpage}{
    accessdate={2021-8-5},
    author={Lagarias, Jeffrey},
    title={Euler's Constant},
    subtitle={Euler's Work and Modern Developments},
    url={https://www.ams.org/journals/bull/2013-50-04/S0273-0979-2013-01423-X/S0273-0979-2013-01423-X.pdf},
    note={Subsection 2.6}
}

\bib{error}{webpage}{
    accessdate={2021-8-5},
    author={Weisstein, Eric},
    title={Erf},
    note={From MathWorld--A Wolfram Web Resource.},
    url={https://mathworld.wolfram.com/Erf.html},
}

\bib{fresnel}{webpage}{
    accessdate={2021-8-5},
    author={Weisstein, Eric},
    title={Fresnel Integrals},
    note={From MathWorld--A Wolfram Web Resource.},
    url={https://mathworld.wolfram.com/FresnelIntegrals.html},
}

\bib{wolframalpha}{webpage}{
    accessdate={2021-8-28},
    author={Wolfram Alpha LLC.},
    url={https://www.wolframalpha.com/input/?i=integrate+\%28x\%5E5\%29\%2F\%28x\%5E7\%2B1\%29+dx+from+x\%3D0.5+to+1.5},
}



\end{biblist}
\end{bibdiv}
\end{document}